\documentclass{amsart}
\usepackage{geometry} 
\usepackage[parfill]{parskip}    
\usepackage{graphicx}
\usepackage{amssymb}
\newcommand{\R}{\mathbb{R}}

\newtheorem{lma}{Lemma}

\newtheorem{thm}{Theorem}

\usepackage{amsfonts}
\begin{document}

\title{A Numerical Minimization Scheme for the Complex Helmholtz Equation}

\author{Russell B. Richins}
\address{Department of Mathematics, University of Utah, Salt Lake City, Utah, 84112-0090}
\email{richins@math.utah.edu}

\author{David C. Dobson}
\address{Department of Mathematics, University of Utah, Salt Lake City, Utah, 84112-0090}
\email{dobson@math.utah.edu}

\subjclass[2000]{Primary 65N30; Secondary 35A15}

\maketitle

\begin{abstract}
We use the work of Milton, Seppecher, and Bouchitt\'{e} on variational principles for waves in lossy media to formulate a finite element method for solving the complex Helmholtz equation that is based entirely on minimization.  In particular, this method results in a finite element matrix that is symmetric positive-definite and therefore simple iterative descent methods and preconditioning can be used to solve the resulting system of equations.  We also derive an error bound for the method and illustrate the method with numerical experiments.
\end{abstract}

\section{Introduction}

Many systems that result in steady-state oscillations can be modeled with the Helmholtz equation, but of particular interest are acoustic waves and transverse electric or transverse magnetic electromagnetic waves in inhomogeneous media.  In each of these situations, the equation of interest can be expressed as  
\begin{equation}\label{acoustic}-\nabla\cdot\rho^{-1}\nabla P-\frac{\omega^2}{\kappa}P=0\end{equation}
for appropriate choices of the complex-valued, spatially dependent material parameters $\rho$ and $\kappa,$ where $\omega > 0$ is the frequency.  The classical methods of deriving a weak form for this equation (with Dirichlet boundary conditions, for example) result in the variational equation
\[\int_\Gamma \left[\rho^{-1}\nabla P\cdot\nabla \overline u-\frac{\omega^2}{\kappa}P\overline u\right] dx=0, \ \ \ \forall u\in H_0^1(\Gamma),\]
which corresponds to a stationary principle, but not a minimization principle.  In \cite{Graeme}, Milton, Seppecher, and Bouchitt\'{e} expand upon the work of Cherkaev and Gibiansky \cite{Cherk} for the conductivity equation to derive variational principles for (\ref{acoustic}) (as a special case of the more general equations of elasticity and electromagnetism) that are true minimization principles,
provided the media are lossy.  The minimization functional corresponds physically to dissipated energy in the system, and is valid 
even for arbitrarily small coefficients of loss. 
While the framework presented in \cite{Graeme} results in nonstandard boundary conditions, Milton and Willis extend the principles to handle the classical Dirichlet and Neumann boundary conditions in \cite{Willis}.   

In this paper we apply the finite element method to these minimization principles and thereby develop a numerical algorithm for solving (\ref{acoustic}) that can take advantage of the many efficient methods available for solving a symmetric, positive-definite system of linear equations.  The outline of the paper is as follows.  Sections~\ref{model} and \ref{bcs} review the general variational formulation and boundary conditions introduced by Milton, Seppecher, Bouchitt\'{e}, and Willis.  In Section~\ref{error_section}, we derive an error bound on certain finite element discretizations of the variational principle.  In Sections~\ref{el_section} and \ref{implement}, we describe a straightforward implementation of the finite element method on a square domain, with Dirichlet boundary conditions.  
In Section~\ref{conditioning_section}, we suggest a preconditioner for solving the resulting symmetric positive definite linear system via the preconditioned conjugate gradient method, and find conditions on the material coefficients under which we expect the best conditioning.  Section~\ref{numerical_section} describes the results of some simple numerical experiments, and illustrates numerical convergence consistent with the error bounds from Section~\ref{error_section}.  Finally, in Section~\ref{robin}, we extend the method to handle Robin boundary conditions, and present some associated numerical examples.

\section{Variational Formulation}

Our model problem is
\begin{equation}\label{model}
\left\{\begin{array}{ll}
\displaystyle -\nabla\cdot\rho^{-1}\nabla P-\frac{\omega^2}{\kappa}P=0 & \mbox{in }\Gamma \\
P=f & \mbox{on }\partial \Gamma \\
\end{array}\right.
\end{equation}
where $\Gamma$ is an open, bounded subset of $\R^d$ ($d=2$ or $3$) with smooth boundary.  For acoustic waves, $\rho$ is the density, $\kappa$ is the bulk modulus, $\omega$ is the frequency, and $P$ is the pressure.  Here $\rho,\ \kappa,$ and $P$ are all complex.  In this section, we focus on Dirichlet boundary conditions for simplicity; Neumann conditions can be handled similarly.  In \cite{Graeme}, it is shown in detail how this and other problems can be formulated as a minimization.  What follows is a brief outline of the general framework.

Let $\mathcal{F}(x)$ and $\mathcal{G}(x)$ be complex-valued fields of the form
\[\mathcal{F}=\left(\begin{array}{c}
F \\
f \\
\end{array}\right)\ \mathcal{G}=\left(\begin{array}{c}
G \\
g \\
\end{array}\right),\]
where $F,\ G:\Gamma\rightarrow\mathbb{C}^d$ and $f,\ g:\Gamma\rightarrow\mathbb{C}$.  Suppose there exists a complex-valued potential $u$ such that 
\[\mathcal{F}=\sqcap u:=\left(\begin{array}{c}
\nabla u \\
u \\
\end{array}\right)\]
and that $\mathcal{G}$ satisfies 
\[h+\sqcup\mathcal{G}=0,\]
where $\sqcup\mathcal{G}:=-\nabla\cdot G+g$ and $h$ is a source term.  Suppose also that $\mathcal{F}$ and $\mathcal{G}$ satisfy the constitutive relation
\begin{equation}\label{constitutive}\mathcal{G}(x)=Z(x)\mathcal{F}(x)\end{equation}
where $Z$ has the form
\[Z=\left(\begin{array}{cc}
L & K \\
K^{T} & M \\
\end{array}\right).\]
Then the constitutive relation along with the differential constraints imply
\begin{equation}\label{one} h+\sqcup(Z\sqcap u)=0\mbox{ or }\nabla\cdot(L\nabla u+Ku)=h+K^{T}\nabla u+Mu\end{equation}
Let $n$ be the unit outward normal on $\partial\Gamma$.  It is shown in \cite{Graeme} that if we are given $u_0'$ and $G_0'\cdot n$ and we specify 
\begin{equation}\label{two}u'=u_0'\mbox{ and } G'\cdot n=G_0'\cdot n\mbox{ on }\partial\Gamma\end{equation}
(herein $'$ denotes the real part of a complex quantity and $''$ the imaginary part), then the solution to (\ref{one}) satisfying the boundary conditions (\ref{two}) is a minimizer of the functional
\[Y(u',G')=\int_\Gamma\left[\left(\begin{array}{c}
\sqcap u' \\
-\mathcal{G}' \\
\end{array}\right)\cdot\mathcal{L}\left(\begin{array}{c}
\sqcap u' \\
-\mathcal{G}' \\
\end{array}\right)+2h''u'\right] dx,\]
where 
\[\mathcal{G}'=\left(\begin{array}{c}
G' \\
\nabla\cdot G'-h \\
\end{array}\right)\mbox{ and }\mathcal{L}=\left(\begin{array}{cc}
Z''+Z'(Z'')^{-1}Z' & Z'(Z'')^{-1} \\
(Z'')^{-1}Z' & (Z'')^{-1} \\
\end{array}\right),\]
provided that $\mathcal{L}$ is positive definite.  An inspection of the constitutive relation shows that $\mathcal{L}$ is positive definite as long as $Z''$ is.  Explicitly, following \cite{Graeme} we see that if we let $\mathcal{F}'$ and $\mathcal{G}'$ be arbitrary, and define $\mathcal{G}''$ and $\mathcal{F}''$ by
\[\left(\begin{array}{c}
\mathcal{G}'' \\
\mathcal{F}'' \\
\end{array}\right)=\mathcal{L}\left(\begin{array}{c}
\mathcal{F}' \\
-\mathcal{G}' \\
\end{array}\right),\]
which is equivalent to 
 \[ \begin{array}{l}
\mathcal{G}'=Z'\mathcal{F}'-Z''\mathcal{F}'' \\
\mathcal{G}''=Z'\mathcal{F}''+Z''\mathcal{F}' \\
\end{array},\]
then
\[\left(\begin{array}{c}
\mathcal{F}' \\
-\mathcal{G}' \\
\end{array}\right)\cdot\mathcal{L}\left(\begin{array}{c}
\mathcal{F}' \\
-\mathcal{G}' \\
\end{array}\right)=\mathcal{F}'\cdot\mathcal{G}''-\mathcal{F}''\cdot\mathcal{G}'\]
\[=\mathcal{F}'\cdot(Z'\mathcal{F}''+Z''\mathcal{F}')-\mathcal{F}''\cdot(Z'\mathcal{F}'-Z''\mathcal{F}'' )\]
\[=\mathcal{F}'\cdot Z''\mathcal{F}'+\mathcal{F}''\cdot Z''\mathcal{F}''.\]
Therefore, $\mathcal{L}$ is positive definite as long as $Z''$ is. 

\section{Boundary Conditions}\label{bcs}

In addition to the conditions 
\[u'=u_0'\mbox{ and }G'\cdot n=G_0'\cdot n\mbox{ on }\partial\Gamma\]
we can also solve the problem for $u'$ and $G'$ with the boundary conditions
\begin{equation}\label{three}u''=u_0''\mbox{ and }G''\cdot n=G_0''\cdot n\mbox{ on }\partial\Gamma,\end{equation}
\[u'=u_0'\mbox{ and }u''=u_0''\mbox{ on }\partial\Gamma,\]
\[\mbox{or }G'\cdot n=G_0'\cdot n\mbox{ and }G''\cdot n=G_0''\cdot n\mbox{ on }\partial\Gamma.\]
The correct variational principles for the last two sets of boundary conditions can be deduced from the formulations for the first two.  The second boundary condition above is a condition on the dual (imaginary) variables $u''$ and $G''$, and therefore it may be enforced through boundary integrals, as follows.

For simplicity, suppose $h=0$.  Let $s\in H^1(\Gamma)$ and $T\in H(\mbox{div},\Gamma)$.  If $u$ and $G$ are such that the differential constraints and constitutive relation are satisfied, then multiplying by $s$ and integrating, we get
\[0=\int_\Gamma(\sqcup \mathcal{G})''s\ dx=\int_\Gamma\left[(-\nabla\cdot G''+g'')s\right]\ dx=\int_\Gamma\left[(-\nabla\cdot G''+g'')s-T\cdot(\nabla u''-\nabla u'')\right]\ dx.\]
Integrating by parts, we find
\[\int_\Gamma\left[G''\cdot\nabla s+g''s-T\cdot\nabla u''-u''\nabla\cdot T\right]\ dx=\int_{\partial\Gamma}\left[sG''\cdot n-u''T\cdot n\right]\ dS.\]
Let $\mathcal{T}=(T,\nabla\cdot T)^T$.  The left-hand side above can be re-written as
\[\int_\Gamma\left(\begin{array}{c}
\mathcal{G}'' \\
\mathcal{F}'' \\
\end{array}\right)\cdot\left(\begin{array}{c}
\sqcap s \\
\mathcal{T} \\
\end{array}\right)\ dx=\int_\Gamma\left(\begin{array}{c}
\sqcap u' \\
\mathcal{G}' \\
\end{array}\right)\cdot\mathcal{L}\left(\begin{array}{c}
\sqcap s \\
\mathcal{T} \\
\end{array}\right)\ dx.\]
The corresponding functional for the boundary condition
\[u''=u_0''\mbox{ and }G''\cdot n=G_0''\cdot n\mbox{ on }\partial\Gamma\]
is then
\[\tilde Y(u',G')=\int_\Gamma\left(\begin{array}{c}
\sqcap u' \\
\mathcal{G}' \\
\end{array}\right)\cdot\mathcal{L}\left(\begin{array}{c}
\sqcap u' \\
\mathcal{G}' \\
\end{array}\right)\ dx+2\int_{\partial\Gamma}\left[u_0''G'\cdot n-u'G_0''\cdot n\right]\ dS.\]
To solve the PDE with the Dirichlet boundary conditions we minimize the functional
\[\hat Y(u',G')=Y(u'+u_0',G')+2\int_{\partial\Gamma}u_0''G'\cdot n\ dS\]
over $u'\in H_0^1(\Gamma)$ and $G'\in H(\mbox{div},\Gamma)$.  To solve the PDE with the Neumann boundary conditions we minimize the functional
\[\check Y(u',G')=Y(u',G'+G_0')-2\int_{\partial\Gamma}u'G_0''\cdot n\ dS\]
over $u'\in H^1(\Gamma)$ and $G'\in H_0(\mbox{div},\Gamma)=\{v\in H(\mbox{div},\Gamma):\langle v\cdot n,w\rangle=0\ \forall\ w\in H_0^1(\Gamma)\}$ (see \cite{Brezzi-Fortin}).

\section{Error Bound}\label{error_section}

In this section we give a bound on the error incurred by solving any of the minimization problems above over a finite dimensional subspace of the specified Sobolev spaces.  We will give a more detailed account of exactly what the finite dimensional space looks like later on; in this section all that will matter is the highest degree of polynomials that the finite dimensional space contains.  We will use the Bramble-Hilbert lemma to give a bound on the error.

Here we will drop the primes used to denote real and imaginary parts.  Note that what follows applies to any of the boundary value problems discussed previously, since the bounds depend only on the corresponding bilinear form.  Throughout this section, $C$ is a constant independent of the solution $(P,v)$ and the grid spacing $h$.

\subsection{Bilinear Form}
Define the bilinear form $B$ by
\begin{equation}\label{bilinear}B(P,v;s,T)=\int_\Gamma\left(\begin{array}{c}
\mathcal{F} \\
-\mathcal{G} \\
\end{array}\right)\cdot\mathcal{L}\left(\begin{array}{c}
\mathcal{S} \\
-\mathcal{T} \\
\end{array}\right)\ dx,\end{equation}
Where, as before, $\mathcal{F}=\sqcap u$, $\mathcal{G}=(G,\nabla\cdot G)^T$, and $\mathcal{S}$ and $\mathcal{T}$ are generated from test function $s\in H^1(\Gamma)$ and $T\in H(\mbox{div},\Gamma)$ in the same fashion.  Assume that there exist constants $\gamma_1,\gamma_2>0$ such that $\mathcal{L}>\gamma_2 I$ and that $[\mathcal{L}(x)]_{ij}\le \gamma_1$ for a.e. $x\in \Gamma$.  Let $V=H_0^1(\Gamma)\times H(\mbox{div},\Gamma)$, endowed with the norm
\[\|(u,G)\|_V=(\|u\|_{H^1(\Gamma)}^2+\|G\|_{H(\mbox{div},\Gamma)}^2)^{\frac{1}{2}}.\]
Then it follows immediately from (\ref{bilinear}) that
\begin{equation}\label{continuous}B(u,G;s,T)\le C\gamma_1\|(u,G)\|_V\|(s,T)\|_V\end{equation}
and
\begin{equation}\label{coercive}B(u,G;u,G)\ge\gamma_2\|(u,G)\|_V^2.\end{equation}

\subsection{Minimization Inequality}
Define an energy by
\[f(s,T)=\frac{1}{2}B(s,T;s,T)-F(s,T),\]
where $F:H^1(\Gamma)\times H(\mbox{div},\Gamma)\rightarrow\R$ (in practice, $F$ is usually composed of terms resulting from an inhomogeneous term and enforcement of the desired boundary conditions). If $(u,G)$ is the minimizer of the energy, then this pair must satisfy the Euler-Lagrange equation
\[B(u,G;s,T)=F(s,T)\ \ \ \forall s\in H_0^1(\Gamma),\ \forall\ T\in H(\mbox{div},\Gamma),\]
so that 
\[f(s,T)=f(u ,G)+\frac{1}{2}B(u-s, G-T;u-s,G-T)\ \ \ \forall\ s\in H_0^1(\Gamma)\ \forall\ T\in H(\mbox{div},\Gamma).\]
Consider a finite dimensional subspace $V_N=V_{N1}\times V_{N2}$ of $V$, where $V_{N1}$ is a finite dimensional subspace of $H^1(\Gamma)$ and $V_{N2}$ is a finite dimensional subspace of $H(\mbox{div},\Gamma)$.  If $(u_N,G_N)$ is such that 
\[f(u_N,G_N)=\min_{(s,T)\in V_N} f(s,T),\]
then 
\[\left[B(u-u_N,G-G_N;u-u_N,G-G_N)\right]^{\frac{1}{2}}=\min_{(s,T)\in V_N}\left[B(u-s,G-T;u-s,G-T)\right]^{\frac{1}{2}}.\]
Inequalities (\ref{continuous}) and (\ref{coercive}) imply that 
\[\sqrt{\gamma_2}\|(s,T)\|_V\le\sqrt{B(s,T;s,T)}\le C\sqrt{\gamma_1}\|(s,T)\|_V\ \ \ \forall\ (s,T)\in V,\]
so we have
\begin{equation}\label{min_test}\sqrt{\gamma_2}\|(u,G)-(u_N,G_N)\|_V\le\min_{(s,T)\in V_N} C \sqrt{\gamma_1}\|(u,G)-(s,T)\|_V.\end{equation}
Let $F_1$ be the orthogonal projection from $H^1(\Gamma)$ onto $V_{N1}$.  Since $F_1$ is an orthogonal projection, it has $\|F_1\|_{B(H^1(\Gamma),H^1(\Gamma))}=1$, where $B(H^1(\Gamma),H^1(\Gamma))$ is the set of bounded linear functions from $H^1(\Gamma)$ to $H^1(\Gamma)$.  Also, define an operator $F_2:H(\mbox{div},\Gamma)\rightarrow V_{N2}$ by the solution of the variational inequality
\[\left<F_2G,Q-F_2G\right>_{L^2(\Gamma,\R^d)}\ge\left<G,Q-F_2G\right>_{L^2(\Gamma,\R^d)}\ \forall\ Q\in E_G,\]
over the set $E_G=\{v\in V_{N2}:\|\nabla\cdot v\|_{L^2(\Gamma)}\le \|\nabla\cdot G\|_{L^2(\Gamma)}\}$, which is a closed, convex subset of $L^2(\Gamma,\R^d)$.  We then have 
\[\|F_2 G\|_{L^2(\Gamma,\R^d)}^2=\left<F_2 G,F_2 G\right>_{L^2(\Gamma,\R^d)}\le\left<G,F_2 G\right>_{L^2(\Gamma,\R^d)}\le\|G\|_{L^2(\Gamma,\R^d)}\|F_2 G\|_{L^2(\Gamma,\R^d)}.\]  
If we take $s=F_1u$ and $T=F_2 G$ in (\ref{min_test}), then we have 
\begin{equation}\label{minimization}\|(u,G)-(u_N,G_N)\|_V\le C\|(u-F_1u,G-F_2G)\|_V.\end{equation}

\subsubsection{Seminorm bounds}

We will discretize the domain $\Gamma$ by by subdividing it into smaller regions, each of which can be seen as a suitable shifting and scaling of a reference element.  More precisely, if $\hat e$ is our reference element, there exist affine changes of variables $F_l(x)=Bx+x_l$ such that $F_l(\hat e)=e_l$, where $e_l$ is the $l$th element (subdivision) in the finite element decomposition of $\Gamma$.  In the case of rectangular elements in $\R^d$, for example, we can take $\hat e=(0,1)^d$, and then we have $B=hI_d$.  In this section a hat will denote the corresponding function defined over the reference element.

Let
\begin{equation}\label{seminorm}[u,w]_s=\sum_{|\alpha|=s}\int_{\hat e}D^{\alpha}u\cdot D^{\alpha}w\ dx\ \mbox{ and }\ |w|_s^2=[w,w]_s, \\
\end{equation}
where for vector functions we define
\[D^{\alpha}w=\left(\begin{array}{c}
D^{\alpha}w_1 \\
D^{\alpha}w_2 \\
\vdots \\
D^{\alpha}w_d \\
\end{array}\right).\]
 From \cite{Brezzi-Fortin} we get the inequalities
\begin{equation}\label{inequalities}\begin{array}{ll}
c^{-1}h^{s-\frac{d}{2}}|w|_{s,e_l}\le|\hat w|_s\le ch^{s-\frac{d}{2}}|w|_{s,e_l} \\
h^{s+\frac{d}{2}-1}|q|_{s,e_l}\le|\hat q|_s\le h^{s+\frac{d}{2}-1}|q|_{s,e_l} \\
h^{s+\frac{d}{2}}|\nabla\cdot q|_{s,e_l}\le |\nabla\cdot\hat q|_s\le h^{s+\frac{d}{2}}|\nabla\cdot q|_{s,e_l} \\
\end{array}\end{equation}
for $B=h$, scalar functions $w$, and vector functions $q$, where $w=\hat w\circ F^{-1}$ and $q=\hat q\circ F^{-1}$ and $|\cdot|_{s,e_l}$ denotes (\ref{seminorm}) with $e_l$ in place of $\hat e$.

We now recall the following lemma from \cite{Axelsson-Barker}, which will be used in what follows.
\begin{lma}[Bramble-Hilbert Lemma] For some region $\Omega\subset\R^2$ and some integer $k\ge -1$, let there be given a bounded linear functional
\[f:H^{k+1}(\Omega)\rightarrow\R,\]
satisfying $|f(u)|\le\delta\|u\|_{H^{k+1}(\Omega)}$ for all $u\in H^{k+1}(\Omega)$ for some $\delta$ independent of $u$.  Suppose that $f(u)=0$ for all $u\in P_k(\bar\Omega)$.  Then there exists a constant $C$, dependent only on $\Omega$ such that 
\[|f(u)|\le C\delta|u|_{k+1},\ \ \ u\in H^{k+1}(\Omega).\]
\end{lma}

Let us suppose that $\hat P\in H^{k+1}(\hat e)$ and $\hat v\in H^{j}(\mbox{div},\hat e)=\{q\in H^j(\hat e,\R^d):\nabla\cdot q\in H^j(\hat e)\}$.  For fixed elements $w\in H^s(\hat e)$ and $Q \in H^s(\mbox{div},\hat e)$ define the functionals
\[f_1(u)=[u-F_1u,w]_s,\ f_2(G)=[G-F_2G,Q]_0,\ f_3(\nabla\cdot G)=[\nabla\cdot G-\nabla\cdot F_2G,\nabla\cdot Q]_0,\]
where $s=0$ or $s=1$.  Then, since
\[|f_1(u)|\le |u-F_1u|_s|w|_s\le (|u|_s+|F_1u|_s)|w|_s\le(\|u\|_{H^1(\Gamma)}+\|F_1u\|_{H^1(\Gamma)})|w|_s\]
\[\le 2\|u\|_{H^1(\Gamma)}|w|_s\le2\|u\|_{H^{k+1}(\Gamma)}|w|_s,\]
\[ |f_2(G)|\le |G-F_2G|_0|Q|_0\le (|G|_0+|F_2G|_0)|Q|_0=(\|G\|_{L^2(\Gamma,\R^d)}+\|F_2G\|_{L^2(\Gamma,\R^d)})|Q|_0\]
\[\le2\|G\|_{L^2(\Gamma),\R^d)}|Q|_0\le\|G\|_{H^j(\Gamma,\R^d)}|Q|_0,\]
\[ |f_3(\nabla\cdot\ G)|\le |\nabla\cdot G-\nabla\cdot F_2G|_0 |\nabla \cdot Q|_0\le (|\nabla\cdot G|_0+|\nabla\cdot F_2G|_0)|\nabla\cdot Q|_0\]
\[=(\|\nabla\cdot G\|_{L^2(\Gamma)}+\|\nabla\cdot F_2 G\|_{L^2(\Gamma)})|\nabla\cdot Q|_0\le 2\|\nabla\cdot G\|_{L^2(\Gamma)}|\nabla\cdot Q|_0\le 2\|\nabla\cdot G\|_{H^j(\Gamma)}|\nabla\cdot Q|_0,\]
and $F_1u=u$ for polynomials in $V_{N1}$ and $F_2G=G$ for vectors of polynomials from $V_{N2}$, we can apply the Bramble-Hilbert lemma to find that there exists a constant such that 
\[|f_1(\hat u)|\le C |w|_s |\hat u|_{k+1},\ |f_2(\hat G)|\le C|Q|_0 |\hat G|_{j},\ |f_3(\nabla \cdot \hat G)| \le C |\nabla \cdot Q|_0|\nabla\cdot\hat G|_{j},\]
as long as $k$ and $j$ are small enough so that all polynomials of degree less than or equal to $k$  are contained in the span of the basis functions representing $\hat u$ and all polynomials of degree less than or equal to $j$ are contained in the span of the basis functions representing $\hat G$.  By choosing $w=\hat u-F_1\hat u$ and $Q=\hat G-F_2\hat G$, we find that
\[|\hat u-F_1\hat u|_s \le C |\hat u|_{k+1},\ |\hat G-F_2\hat G|_0\le C |\hat G|_{j},\ |\nabla\cdot\hat G-\nabla\cdot F_2\hat G|_0\le C|\nabla\cdot \hat G|_{j}.\]

Employing (\ref{inequalities}), we find that for $h\le 1$,
\[|u-F_1u|_{s,e_l}\le C h^{\frac{d}{2}-s}|\hat u-F_1\hat u|_s\le C h^{\frac{d}{2}-s}|\hat u|_{k+1}\le C h^{k-s+1}|u|_{k+1,e_l},\]
\[|G-F_2G|_{0,e_l}\le h^{1-\frac{d}{2}}|\hat G-F_2\hat G|_0\le h^{1-\frac{d}{2}}C|\hat G|_{j}\le Ch^{j}|G|_{j,e_l},\]
\[|\nabla\cdot G-\nabla\cdot F_2G|_{0,e_l}\le h^{-\frac{d}{2}}|\nabla \cdot \hat G-\nabla\cdot F_2\hat G|_{0} \le h^{-\frac{d}{2}}C|\nabla\cdot\hat G|_{j}\le C h^{j}|\nabla \cdot G|_{j,e_l}.\]
Returning to inequality (\ref{minimization}), we have 
\[\|(u,G)-(u_N,G_N)\|_V^2\le C\|(u,G)-(F_1u,F_2G)\|_V^2\]
\[=C \sum_l\left[|u-F_1 u|_{0,e_l}^2+|u-F_1u|_{1,e_l}^2+|G-F_2 G|_{0,e_l}^2+|\nabla\cdot v-\nabla\cdot F_2 G|_{0,e_l}^2\right]\] 
\[\le C \sum_l \left[h^{2k+2}|u|_{k+1,e_l}^2+h^{2k}|u|_{k+1,e_l}^2+h^{2j}|G|_{j,e_l}^2+h^{2j}|\nabla \cdot G|_{j,e_l}^2\right] \]
\[ \le C (h^{2k} |u|_{k+1,\Gamma}^2+h^{2j}(|G|_{j,\Gamma}^2+|\nabla\cdot G|_{j,\Gamma}^2)).\]

Let $P_k(\bar \Gamma)$ denote all polynomials of degree less than or equal to $k$ on $\bar \Gamma$.  We have now proved
\begin{thm}\label{result}If the solution $(u,G)\in H^{k+1}(\Gamma)\times H^{j+1}(\mbox{div},\Gamma)$ and the finite element subspace used in the numerical method contains $P_k(\bar \Gamma) \times P_j(\bar \Gamma)) \times P_j(\bar \Gamma)$, then there is a constant $C$ such that the error satisfies
\[\|(u,G)-(u_N,G_N)\|_V^2\le C (h^{2k} |u|_{k+1,\Gamma}^2+h^{2j}(|G|_{j,\Gamma}^2+|\nabla\cdot G|_{j,\Gamma}^2)),\]
where $h\le 1$ is the grid spacing.
\end{thm}

\subsection{Regularity}
In order for the error bound to be meaningful, we must have $k,j \ge 1$ in Theorem~\ref{result}, which means that at least
\[u\in H^2(\Gamma)\mbox{ and }G\in H^1(\mbox{div},\Gamma).\]
In the notation of the acoustic equation, if $\rho^{-1}$ is positive definite, bounded, and $C^1$, then classical elliptic regularity theory such as in \cite{Evans} guarantees that $P'\in H^2(\Gamma)$.  Also since 
\[v=-\frac{i}{\omega}\rho^{-1}\nabla P,\]
we have that $v\in (H^1(\Gamma))^2$, and multiplying the acoustic equation through by $-1/\omega$ tells us that 
\[\nabla\cdot v=\frac{i\omega}{\kappa}P,\]
so $\nabla\cdot v\in H^1(\Gamma)$ as long as $\kappa$ is at least $C^1$.

It would be more satisfying (and useful in other contexts) to have a regularity theory derived from the weak form of the equations presented herein, and this is a current topic of inquiry for the authors.

\section{Euler-Lagrange Equation for the Model Problem}\label{el_section}
For our model in the development of the numerical method, we will focus on the Dirichlet problem with functional
\[\hat Y(u',G')=\int_\Gamma \left(\begin{array}{c}
\sqcap u' \\
-\mathcal{G}' \\
\end{array}\right)\cdot\mathcal{L}\left(\begin{array}{c}
\sqcap u' \\
-\mathcal{G}' \\
\end{array}\right)\ dx+2\int_{\partial\Gamma}G'\cdot n u''\ dS.\]
Suppose that $u'$ and $u''$ satisfy (\ref{three}) and $(u',G')$ minimizes $\hat Y$ over all $u'\in u_0'+H_0^1(\Gamma)$ and $G'\in H(\mbox{div},\Gamma)$.  Then if we take any functions $s\in H_0^1(\Gamma)$ and $T\in H(\mbox{div},\Gamma)$ and let $\mathcal{T}=(T,\nabla\cdot T)^T$, we have that
\[\hat Y(u'+ts,G'+tT)=\int_\Gamma\left(\begin{array}{c}
\sqcap u'+t\sqcap s \\
-\mathcal{G}'-t\mathcal{T} \\
\end{array}\right)\cdot\mathcal{L}\left(\begin{array}{c}
\sqcap u'+t\sqcap s \\
-\mathcal{G}'-t\mathcal{T} \\
\end{array}\right) dx\]
\[+2\int_{\partial\Gamma}(G'+tT)\cdot n u''\ dS \]
has a minimum at $t=0$.  Therefore,
\begin{equation}\label{euler-lagrange}0=2\int_\Gamma\left(\begin{array}{c}
\sqcap u' \\
-\mathcal{G}' \\
\end{array}\right)\cdot\mathcal{L}\left(\begin{array}{c}
\sqcap s \\
-\mathcal{T} \\
\end{array}\right)\ dx+2\int_{\partial\Gamma}T\cdot n u''\ dS.\end{equation}
This is the weak form of the equation that we want to solve for $u'$.  In the case of the acoustic equation (\ref{model}), we have 
\[u=P,\ L=-\rho^{-1},\ K=0,\  M=\omega^2/\kappa,\ h=0,\ v=(-i/\omega)\rho^{-1}\nabla P,\ G=-i\omega v,\]
so we can rewrite
\begin{align}\label{yhat}\hat Y(P',v'')&=\int_\Gamma\left[\left(\begin{array}{c}
\nabla P' \\
-\omega v'' \\
\end{array}\right)\cdot\mathcal{R}\left(\begin{array}{c}
\nabla P' \\
-\omega v'' \\
\end{array}\right)+\left(\begin{array}{c}
\omega P' \\
-\nabla\cdot v'' \\
\end{array}\right)\cdot\mathcal{K}\left(\begin{array}{c}
\omega P' \\
-\nabla\cdot v''\\
\end{array}\right)\right] dx\notag \\
&+2\int_{\partial\Gamma}\omega v''\cdot nP''\ dS,\end{align}
where $r=-\rho^{-1}$, $k=\kappa^{-1}$, and
\[\mathcal{R}=\left(\begin{array}{cc}
r''+r'(r'')^{-1}r' & r'(r'')^{-1} \\
(r'')^{-1}r' & (r'')^{-1} \\
\end{array}\right),\ \mathcal{K}=\left(\begin{array}{cc}
k''+(k')^2/k'' & k'/k'' \\
k'/k'' & 1/k'' \\
\end{array}\right).\]
The requirement that $Z''$ be positive definite translates to the requirement that 
\begin{equation}\label{coerc_bounds}\rho''>\alpha I,\ \kappa''<-\beta,\ \alpha,\beta > 0.\end{equation}

Making the substitutions in (\ref{euler-lagrange}) for the acoustic equation, we find that the Euler-Lagrange equation becomes
\begin{equation}\label{four}0=\int_\Gamma\left[\left(\begin{array}{c}
\nabla P' \\
-\omega v'' \\
\end{array}\right)\cdot\mathcal{R}\left(\begin{array}{c}
\nabla s \\
-\omega T \\
\end{array}\right)+\left(\begin{array}{c}
\omega P' \\
-\nabla\cdot v'' \\
\end{array}\right)\cdot\mathcal{K}\left(\begin{array}{c}
\omega s \\
-\nabla\cdot T \\
\end{array}\right)\right] dx+\int_{\partial\Gamma}\omega T\cdot nP''\ dS\end{equation}
for any $s\in H_0^1(\Gamma)$ and any $T\in H(\mbox{div},\Gamma)$.

\section{Implementation of the Finite Element Method}\label{implement}

Our goal is to test the efficacy of this new variational principle, using a simple, explicit finite element implementation.  Let us assume that $d=2$ and $\Gamma=(0,1)^2$.  In order to find a numerical solution for $P'$, we introduce an $N\times N$ computational grid  with equally spaced nodes $(x_j,y_t)$ for $t,j=1,2,\ldots,N$ and grid spacing $h=1/(N-1)$.  We also introduce the finite element spaces
\[\begin{array}{l}
\displaystyle \Psi=\mbox{span}\left\{\left(1-\frac{|x-x_j|}{h}\right)\left(1-\frac{|y-y_t|}{h}\right)\chi_{tj}:2\le t,j\le N-1\right\} \\
 \Phi_1=\mbox{span}\left\{\left(\begin{array}{c}
\displaystyle \left(1-\frac{|x-x_j|}{h}\right)\left(1-\frac{|y-y_t|}{h}\right) \\
0 \\
\end{array}\right)\chi_{tj}:1\le t,j\le N\right\} \\
\Phi_2=\mbox{span}\left\{\left(\begin{array}{c}
0 \\
\displaystyle \left(1-\frac{|x-x_j|}{h}\right)\left(1-\frac{|y-y_t|}{h}\right) \\
\end{array}\right)\chi_{tj}:1\le t,j\le N\right\} \\
\end{array},\]
where 
\[\chi_{tj}(x,y)=\left\{\begin{array}{ll}
1 & \mbox{if } |x-x_j|,|y-y_t|\le h \\
0 & \mbox{otherwise} \\
\end{array}\right. .\]
The bases of each of the finite element spaces are built from simple piecewise bilinear elements.

We can re-index these elements with a single index by setting
\[\begin{array}{l}
\displaystyle \psi_k=\left(1-\frac{|x-x_j|}{h}\right)\left(1-\frac{|y-y_t|}{h}\right)\chi_{tj},\mbox{ where }k=(t-2)(N-2)+j-1,\ k=1,\ldots,(N-2)^2, \\
\phi_{1k}=\left(\begin{array}{c}
\displaystyle \left(1-\frac{|x-x_j|}{h}\right)\left(1-\frac{|y-y_t|}{h}\right) \\
0 \\
\end{array}\right)\chi_{tj}\mbox{ where }k=(t-1)N+j,\ k=1,\ldots,N^2, \\
\phi_{2k}=\left(\begin{array}{c}
0 \\
\displaystyle \left(1-\frac{|x-x_j|}{h}\right)\left(1-\frac{|y-y_t|}{h}\right) \\
\end{array}\right)\chi_{tj}\mbox{ where }k=(t-1)N+j,\ k=1,\ldots,N^2. \\
\end{array}\]
We assume that our finite element solution has the form
\[\left(\begin{array}{c}
P' \\
v'' \\
\end{array}\right)=\left(\begin{array}{c}
\psi_R+\sum_{k=1}^{(N-2)^2}\delta_k\psi_k \\
\sum_{k=1}^{N(N-1)}\beta_k\phi_{1k}+\sum_{k=1}^{N(N-1)}\gamma_k\phi_{2k} \\
\end{array}\right).\]
Here $\psi_R$ is any function that satisfies the desired Dirichlet boundary condition for $P'$.  Making this substitution into (\ref{four}), we get
\begin{flalign}\int_\Gamma & \left[\left(\begin{array}{c}
\sum\delta_k\nabla\psi_k \\
-\omega\sum\beta_k\phi_{1k}-\omega\sum\gamma_k\phi_{2k} \\
\end{array}\right)\cdot\mathcal{R}\left(\begin{array}{c}
\nabla s \\
-\omega T \\
\end{array}\right)\right. & \notag\end{flalign}
\[\left.+\left(\begin{array}{c}
\omega\sum\delta_k\psi_k \\
-\sum\beta_k\nabla\cdot\phi_{1k}-\sum\gamma_k\nabla\cdot\phi_{2k} \\
\end{array}\right)\cdot\mathcal{K}\left(\begin{array}{c}
\omega s \\
-\nabla\cdot T \\
\end{array}\right)\right] dx\]
\begin{flalign} &=-\int_\Gamma\left[\left(\begin{array}{c}
\nabla\psi_0 \\
0 \\
\end{array}\right)\cdot\mathcal{R}\left(\begin{array}{c}
\nabla s \\
-\omega T \\
\end{array}\right)+\left(\begin{array}{c}
\omega\psi_0 \\
0 \\
\end{array}\right)\cdot\mathcal{K}\left(\begin{array}{c}
\omega s \\
-\nabla\cdot T \\
\end{array}\right)\right] dx & -\int_\Gamma\left[\omega\nabla\psi_I\cdot T+\omega\psi_I\nabla\cdot T\right] dx,\notag\end{flalign}
where we have used the divergence theorem on the boundary integral, $\psi_I$ is any function on $\Gamma$ satisfying the desired Dirichlet boundary condition for $P''$, and $s\in H_0^1(\Gamma),\ T\in H(\mbox{div},\Gamma)$ are arbitrary.  In particular, this must hold when
\[s=\psi_k,\ T=0\mbox{ for }k=1,\ldots,(N-2)^2\]
\[s=0,\ T=\phi_{1k}\mbox{ for }k=1,\ldots,N(N-1)\]
\[s=0,\ T=\phi_{2k}\mbox{ for }k=1,\ldots,N(N-1).\]
This gives rise to a system of equations of the form $A\alpha=b$, where $A$ has the block form
\begin{equation}\label{block}A=\left(\begin{array}{ccc}
A_1 & A_4 & A_6 \\
A_4 & A_2 & A_5 \\
A_6 & A_5 & A_3 \\
\end{array}\right)\end{equation}
and the blocks have entries
\begin{eqnarray}\label{block_entries}(A_1)_{tj}&=&\displaystyle\int_\Gamma\left[\left(\begin{array}{c}
\nabla \psi_t \\
0 \\
\end{array}\right)\cdot\mathcal{R}\left(\begin{array}{c}
\nabla\psi_j \\
0 \\
\end{array}\right)+\left(\begin{array}{c}
\omega\psi_t \\
0 \\
\end{array}\right)\cdot\mathcal{K}\left(\begin{array}{c}
\omega\psi_j \\
0 \\
\end{array}\right)\right] dx \notag \\
(A_2)_{tj}&=&\displaystyle\int_\Gamma\left[\left(\begin{array}{c}
0 \\
-\omega\phi_{1t} \\
\end{array}\right)\cdot\mathcal{R}\left(\begin{array}{c}
0 \\
-\omega\phi_{ij} \\
\end{array}\right)+\left(\begin{array}{c}
0 \\
-\nabla\cdot\phi_{1t} \\
\end{array}\right)\cdot\mathcal{K}\left(\begin{array}{c}
0 \\
-\nabla\cdot\phi_{1j} \\
\end{array}\right)\right] dx \notag \\
(A_3)_{tj}&=&\displaystyle\int_\Gamma\left[\left(\begin{array}{c}
0 \\
-\omega\phi_{2t} \\
\end{array}\right)\cdot\mathcal{R}\left(\begin{array}{c}
0 \\
-\omega\phi_{2j} \\
\end{array}\right)+\left(\begin{array}{c}
0 \\
-\nabla\cdot\phi_{2t} \\
\end{array}\right)\cdot\mathcal{K}\left(\begin{array}{c}
0 \\
-\nabla\cdot\phi_{2j} \\
\end{array}\right)\right] dx  \\
(A_4)_{tj}&=&\displaystyle\int_\Gamma\left[\left(\begin{array}{c}
0 \\
-\omega\phi_{1t} \\
\end{array}\right)\cdot\mathcal{R}\left(\begin{array}{c}
\nabla\psi_j \\
0 \\
\end{array}\right)+\left(\begin{array}{c}
0 \\
-\nabla\cdot\phi_{1t} \\
\end{array}\right)\cdot\mathcal{K}\left(\begin{array}{c}
\omega\psi_j \\
0 \\
\end{array}\right)\right] dx. \notag \\
(A_5)_{tj}&=&\displaystyle\int_\Gamma\left[\left(\begin{array}{c}
0 \\
-\omega\phi_{2t} \\
\end{array}\right)\cdot\mathcal{R}\left(\begin{array}{c}
0 \\
-\omega\phi_{1j} \\
\end{array}\right)+\left(\begin{array}{c}
0 \\
-\nabla\cdot\phi_{2t} \\
\end{array}\right)\cdot\mathcal{K}\left(\begin{array}{c}
0 \\
-\nabla\cdot\phi_{1j} \\
\end{array}\right)\right] dx \notag \\
(A_6)_{tj}&=&\displaystyle\int_\Gamma\left[\left(\begin{array}{c}
0 \\
-\omega\phi_{2t} \\
\end{array}\right)\cdot\mathcal{R}\left(\begin{array}{c}
\nabla\psi_j \\
0 \\
\end{array}\right)+\left(\begin{array}{c}
0 \\
-\nabla\cdot\phi_{2t} \\
\end{array}\right)\cdot\mathcal{K}\left(\begin{array}{c}
\omega\psi_j \\
0 \\
\end{array}\right)\right] dx. \notag 
\end{eqnarray}
The right-hand side vector $b$ is partitioned as
\[b=\left(\begin{array}{c}
b_1 \\
b_2 \\
b_3 \\
\end{array}\right),\]
where 
\begin{equation}\label{five}\begin{array}{ll}
(b_1)_k= & -\displaystyle\int_\Gamma\left[\left(\begin{array}{c}
\nabla\psi_R \\
0 \\
\end{array}\right)\cdot\mathcal{R}\left(\begin{array}{c}
\nabla\psi_k \\
0 \\
\end{array}\right)+\left(\begin{array}{c}
\omega\psi_R \\
0 \\
\end{array}\right)\cdot\mathcal{K}\left(\begin{array}{c}
\omega\psi_k \\
0 \\
\end{array}\right)\right] dx \\
(b_2)_k= & -\displaystyle\int_\Gamma\left[\left(\begin{array}{c}
\nabla\psi_R \\
0 \\
\end{array}\right)\cdot\mathcal{R}\left(\begin{array}{c}
0 \\
-\omega\phi_{1k} \\
\end{array}\right)+\left(\begin{array}{c}
\omega\psi_R \\
0 \\
\end{array}\right)\cdot\mathcal{K}\left(\begin{array}{c}
0 \\
-\nabla\cdot\phi_{1k} \\
\end{array}\right)\right] dx \\
 & -\displaystyle\int_\Gamma\left[\omega\nabla\psi_I\cdot\phi_{1k}+\omega\psi_I\nabla\cdot\phi_{1k}\right] dx \\
(b_3)_k= & -\displaystyle\int_\Gamma\left[\left(\begin{array}{c}
\nabla\psi_R \\
0 \\
\end{array}\right)\cdot\mathcal{R}\left(\begin{array}{c}
0 \\
-\omega\phi_{2k} \\
\end{array}\right)+\left(\begin{array}{c}
\omega\psi_R \\
0 \\
\end{array}\right)\cdot\mathcal{K}\left(\begin{array}{c}
0 \\
-\nabla\cdot\phi_{2k} \\
\end{array}\right)\right] dx \\
 & -\displaystyle\int_\Gamma\left[\omega\nabla\psi_I\cdot\phi_{2k}+\omega\psi_I\nabla\cdot\phi_{2k}\right] dx. \\
\end{array}\end{equation}

The method for solving for $P'$ and $v''$ can be easily modified to solve for $P''$ and $v'$.  In this case the weak equation is
\[\int_\Gamma\left[\left(\begin{array}{c}
\nabla P'' \\
\omega v' \\
\end{array}\right)\cdot\mathcal{R}\left(\begin{array}{c}
\nabla s \\
\omega T \\
\end{array}\right)+\left(\begin{array}{c}
-\omega P'' \\
-\nabla\cdot v' \\
\end{array}\right)\cdot\mathcal{K}\left(\begin{array}{c}
-\omega s \\
-\nabla\cdot T \\
\end{array}\right)\right] dx+\int_{\partial\Gamma}\omega T\cdot nP'\ dS,\]
and all the methods above still apply.  In fact, to obtain the new matrix for this formulation, we simply change the signs of the blocks $A_4$ and $A_6$, and the changes in $b$ are mostly reversing signs and the roles of the two auxiliary functions $\psi_R$ and $\psi_I$. 

\subsection{Other Discretizations}\label{other_disc}
Along with the discretization described in Section~\ref{implement}, we have experimented with two other implementations in which different basis functions are used to represent the variable $v$.  The first of these uses the Raviart-Thomas $RT_{[0]}$ elements described in \cite{Brezzi-Fortin}.  We found that the resulting finite element matrix is much more poorly scaled, with a condition number approximately twice as large as that obtained with the nodal bilinear basis.  The second method uses the $RT_{[1]}$ elements (also described in \cite{Brezzi-Fortin}).  In this case, the higher-order basis functions obviously result in a somewhat less-sparse finite element matrix, and the condition number is approximately the same as that obtained with the all-bilinear discretization.

\section{Conditioning}\label{conditioning_section}

 \begin{figure}[tbp]
   \centering
   \includegraphics[width=7in]{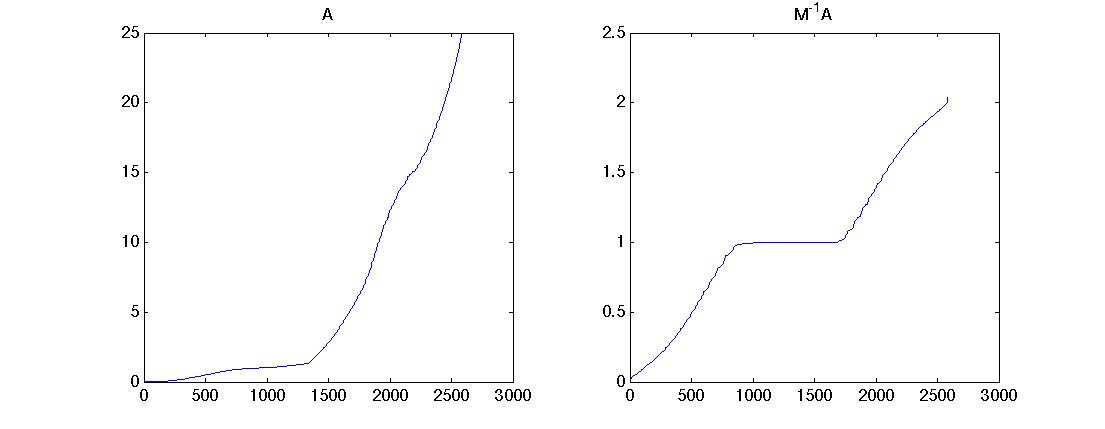} 
   \caption{The distribution of the eigenvalues of $A$ and the real parts of the eigenvalues of $M^{-1}A$ for $N=30$.}
   \label{evals}
\end{figure}

As was mentioned, perhaps the greatest numerical advantage to having a minimization formulation for the Helmholtz equation is that the matrix produced by the finite element method is symmetric positive definite.  This allows for the use of methods such as the conjugate gradient method to solve the system.  Of course, the use of a preconditioning matrix in the conjugate gradient method can speed up the convergence considerably, which is especially important when solving the relatively large sparse systems generated by the finite element approach outlined above.

In our approach, there are three basic types of elements used: bilinear elements, first component bilinear vector elements, and second component bilinear vector elements.  Each of these types of elements interacts with all of the other types, and these interactions are what give rise to the blocks in (\ref{block}).  Assuming that interactions among similar element types are most important, we choose the block Jacobi preconditioner 
\[M=\left(\begin{array}{ccc}
A_1 & 0 & 0 \\
0 & A_2 & 0 \\
0 & 0 & A_3 \\
\end{array}\right).\]
Among all block diagonal preconditioners of this form, this choice of $M$ minimizes the condition number of $M^{-\frac{1}{2}}AM^{-\frac{1}{2}}$ to within a factor of 3 of its minimum \cite{Demmel}.  

\begin{figure}[tbp]
   \centering
   \includegraphics[width=4in]{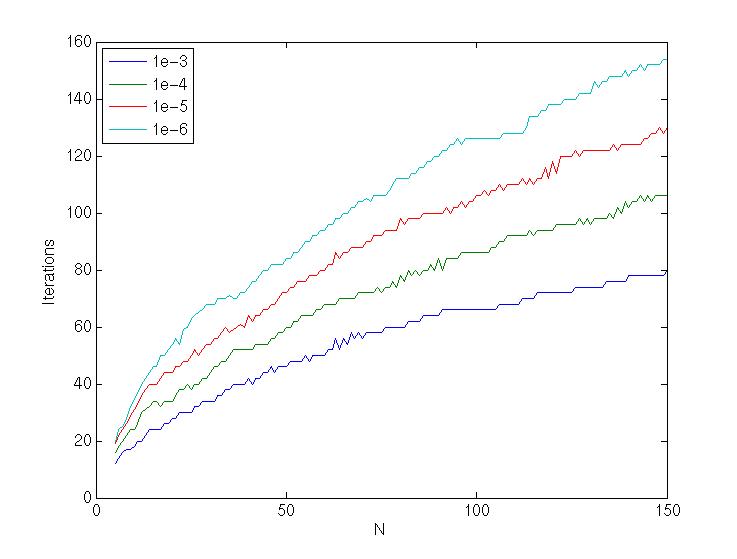} 
   \caption{The growth of the number of PCG iterations required to solve a given problem with grid size for several error tolerances (outer implementation of PCG only).}
   \label{iters}
\end{figure}

As one of the steps in the preconditioned conjugate gradient method (PCG), \cite{Demmel_book}, a system of the form $Mr=y$ must be solved.  In order to make solving this problem more efficient, we precondition the matrix $M$ and use conjugate gradient to solve this system as well.  The preconditioner used in this inner implementation of PCG was an incomplete Choleski factorization of M.  Figure \ref{evals} shows the distribution of the eigenvalues of the matrix $A$ before and after preconditioning for $N=30$.    In Figure \ref{iters}, we see the how the number of PCG iterations grows with $N$ for several error tolerances.

A  key component in ensuring that the system $A\alpha=b$ is well conditioned is for the matrix $\mathcal{L}$ (or equivalently $\mathcal{R}$ and $\mathcal{K}$) to have a coercivity constant that is as large as possible.  For this reason, we expect better numerical results when the eigenvalues of $\mathcal{L}$ are bounded well away from zero. In the case of the Helmholtz equation, the matrix $Z$ is diagonal, say $Z=\mbox{diag}(c_1,\ldots,c_{d+1})$, which makes it possible to calculate the eigenvalues of $\mathcal{L}$.  If $D$ is an invertible matrix, then we may factor a block diagonal matrix  
\[\left(\begin{array}{cc}
A & B \\
C & D \\
\end{array}\right)\]
as
\[\left(\begin{array}{cc}
A & B \\
C & D \\
\end{array}\right)=\left(\begin{array}{cc}
I & B \\
0 & D \\
\end{array}\right)\left(\begin{array}{cc}
A-BD^{-1}C & 0 \\
D^{-1}C & I \\
\end{array}\right),\]
which implies that 
\[\mbox{det}\left(\begin{array}{cc}
A & B \\
C & D \\
\end{array}\right)=\mbox{det}(D)\mbox{det}(A-BD^{-1}C).\]
Therefore, 
\[\mbox{det}(\mathcal{L}-\lambda I)=(-1)^{d+1}\mbox{det}\left(\begin{array}{cc}
(Z'')^{-1}Z' & (Z'')^{-1}-\lambda I \\
Z''+Z'(Z'')^{-1}Z'-\lambda I & Z'(Z'')^{-1} \\
\end{array}\right)\]
\[=(-1)^{d+1}\mbox{det}(Z'(Z'')^{-1})\mbox{det}((Z'')^{-1}Z'+[-(Z')^{-1}+\lambda Z''(Z')^{-1}][Z''+Z'(Z'')^{-1}Z'-\lambda I])\]
\[=(-1)^{d+1}\mbox{det}(Z'(Z'')^{-1})\mbox{det}(\lambda^2[-Z''(Z')^{-1}]+\lambda[(Z')^{-1}+Z''(Z')^{-1}Z''+Z']-(Z')^{-1}Z'').\]
In the case of diagonal $Z$, this implies that 
\[\lambda=\frac{-a_j\pm\sqrt{a_j^2-b_j^2}}{-b_j}\ \ \ j=1,\ldots,d+1,\]
where
\[a_j=\frac{1}{c_j'}+\frac{(c_j'')^2}{c_j'}+c_j'\ \mbox{ and }\ b_j=2\frac{c_j''}{c_j'}.\]
If $Z'=0$, then $\mathcal{L}$ is diagonal, and its eigenvalues are those of $Z''$ and $(Z'')^{-1}$. 

The above analysis tells us that the finite element problem will be better conditioned for those problems where the coefficients $\rho$ and $\kappa$ are such that $Z$ is close to $Ii$, i.e. $\rho=iI$ and $\kappa=-Ii$ (this would correspond to the limiting case where $a_j=b_j$).  In many cases when we are presented with a problem where the coercivity constant for $\mathcal{L}$ is small, we can apply an appropriate rotation and scaling to the problem in order to get a finite element matrix that is better conditioned.  By multiplying the problem (\ref{model}) through by a complex constant $re^{i\theta}$, we effectively replace $Z$ with $re^{i\theta}Z$, so we should choose $r$ and $\theta$ so that $re^{i\theta}Z$ is as close as possible to $iI$.  However, this may not always be possible, for example, when an isotropic $\rho(x)$ oscillates between values in the upper half of the complex plane that are close to $1$ and $-1$.

\section{Numerical Results}\label{numerical_section}

\begin{figure}[t]
   \centering
   \includegraphics[width=2.953in]{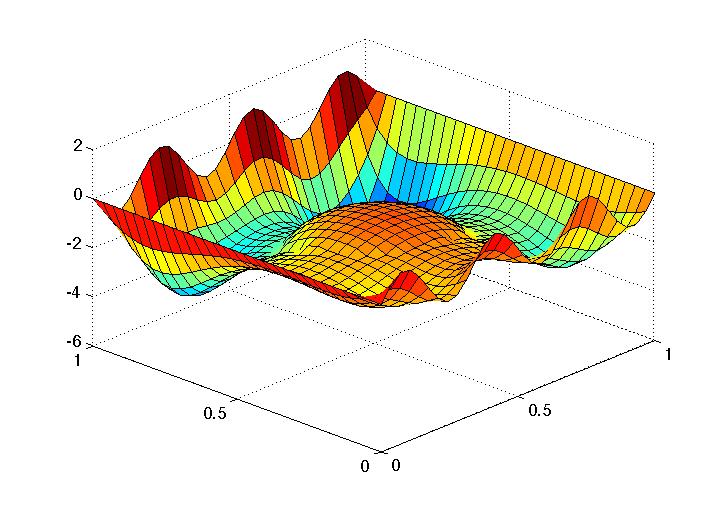} 
   \includegraphics[width=2.953in]{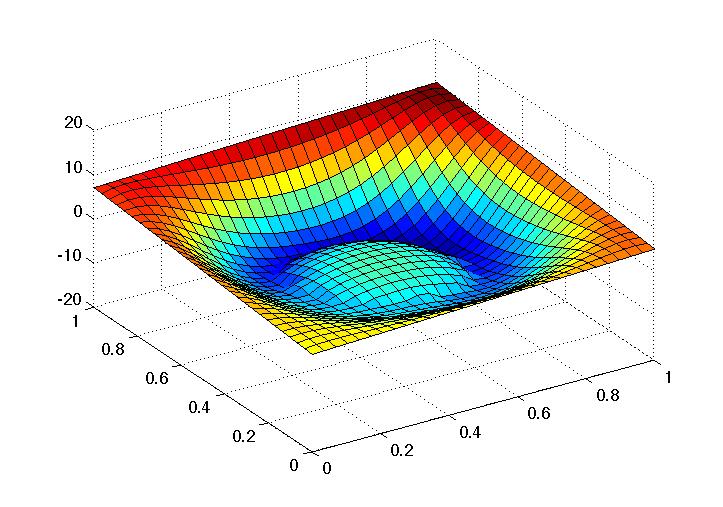} 
   \caption{The solutions $P'$ and $P''$,  with $\omega=10,\psi_R=\sin(6\pi x)\cos(3\pi y),\ \psi_I=3x+5y+2,\ N=30$.  There is a circular inclusion in the center of the domain with  $\rho=.01+.001i$ and $\kappa=.01-.003i$ outside the inclusion and $\rho=-5+5i$ and $\kappa=4-4i$ inside the inclusion.}
     \label{example_fig}
\end{figure}

As an example, we demonstrate the error bound on the problem (\ref{acoustic}), with parameters $\rho=(-5+5i)I,\ \kappa=4-4i$ and $\omega=2$.  A solution is $P(x,y)=e^{2ix-3y}$.  In this example we took 
\[\psi_R=\mbox{Re}(e^{2ix-3y})+\sin(\pi x)\sin(\pi y),\ \psi_I=\mbox{Im}(e^{2ix-3y})+\sin(\pi x)\cdot 3\sin(\pi y))\]
and solved the problem on grids with $N=3,\ldots,100$.  Table 1 shows the error in the finite element solution for various values of $N$.  The errors were calculated using the trapezoidal rule with function evaluations on a grid with size $N=1500$.  Figure~\ref{example_fig} demonstrates the method on a problem with non-constant coefficients, where the dissipation in the material is higher inside a disk centered in the unit square.  The boundary conditions for the real part are oscillatory, while the boundary conditions for the imaginary part are simply an affine function.

\begin{table}
\label{error_table}
\begin{center}
\begin{tabular}{|c|c|c|}
\hline
N & h & $\|(P-P_N,v-v_N)\|_V$ \\
\hline
30 & 0.0345 & $6.6162\times 10^{-4} $\\
40 & 0.0256 & $3.6692\times 10^{-4}$ \\
50 & 0.0204 & $ 2.3252\times 10^{-4}$ \\
60 & 0.0169 & $1.6026\times 10^{-4}$ \\
70 & 0.0145 & $1.1722\times 10^{-4}$ \\
80 & 0.0127 & $8.9706\times 10^{-5}$ \\
90 & 0.0112 & $7.0686\times 10^{-5} $\\
100 & 0.0101 & $5.7037\times 10^{-5}$ \\
\hline
\end{tabular}
\end{center}
\caption{The error in the finite element solution for various values of the grid size $N$.}
\end{table}

\section{Robin Boundary Conditions}\label{robin}

\subsection{Problem Formulation}

Another boundary condition that often appears is the Robin problem
\[\left\{\begin{array}{ll}
\displaystyle -\nabla\cdot\rho^{-1}\nabla P-\frac{\omega^2}{\kappa}P=0 & \mbox{in}\ \Gamma, \\
P+av\cdot n=g & \mbox{on}\ \partial\Gamma, \\
\end{array}\right.\]
where $a\in\mathbb{C}$.  In order to deal with this boundary condition, which concerns both real and imaginary parts of the variables $P$ and $v$ simultaneously, we start with the minimization functional for the natural boundary conditions
\[Y(P',v'')+2\omega\int_{\partial\Gamma}\left[P'v'\cdot n+P''v''\cdot n\right]dS.\]
The Euler-Lagrange Equation for the corresponding variational principle is 
\[B(P',v'',s,T)=-\omega\int_{\partial\Gamma}\left[sv'\cdot n+P''T\cdot n\right] dS.\]
where the bilinear form $B$ is defined in (\ref{bilinear}).  Notice that we can write the surface integral above as
\[-\omega\int_{\partial\Gamma}\left(\begin{array}{c}
v'\cdot n \\
P'' \\
\end{array}\right)\cdot\left(\begin{array}{c}
s \\
T\cdot n \\
\end{array}\right)\ dS.\]
The vector on the right contains the primary variables for which we would like to solve, and the vector on the left contains the dual variables which we would like to eliminate using the Robin boundary condition.  In terms of the vectors above, we can express the Robin condition as 
\[M_1\left(\begin{array}{c}
P' \\
v''\cdot n \\
\end{array}\right)+M_2\left(\begin{array}{c}
v'\cdot n \\
P'' \\
\end{array}\right)=\left(\begin{array}{c}
g' \\
g'' \\
\end{array}\right),\]
where 
\[M_1=\left(\begin{array}{cc}
1 & -a'' \\
0 & a' \\
\end{array}\right)\ \mbox{and}\ M_2=\left(\begin{array}{cc}
a' & 0 \\
a'' & 1 \\
\end{array}\right).\]
Rearranging, we find that
\[\left(\begin{array}{c}
v'\cdot n \\
P'' \\
\end{array}\right)=M_2^{-1}\left(\begin{array}{c}
g' \\
g'' \\
\end{array}\right)-M_2^{-1}M_1\left(\begin{array}{c}
P' \\
v''\cdot n \\
\end{array}\right),\]
so the surface integral term becomes
\[-\omega\int_{\partial\Gamma}\left[M_2^{-1}\left(\begin{array}{c}
g' \\
g'' \\
\end{array}\right)-M_2^{-1}M_1\left(\begin{array}{c}
P' \\
v''\cdot n \\
\end{array}\right)\right]\cdot\left(\begin{array}{c}
S \\
T\cdot n \\
\end{array}\right)\ dS\]
\[=-\omega\int_{\partial\Gamma}M_2^{-1}\left(\begin{array}{c}
g' \\
g'' \\
\end{array}\right)\cdot\left(\begin{array}{c}
P' \\
v''\cdot n \\
\end{array}\right)\ dS+\omega\int_{\partial\Gamma}M_2^{-1}M_1\left(\begin{array}{c}
P' \\
v''\cdot n \\
\end{array}\right)\cdot\left(\begin{array}{c}
P' \\
v''\cdot n \\
\end{array}\right)\ dS.\]
The new Euler-Lagrange equation for the Robin boundary condition is therefore
\[B(P',v'';s,T)-\omega\int_{\partial\Gamma}M_2^{-1}M_1\left(\begin{array}{c}
P' \\
v''\cdot n \\
\end{array}\right)\cdot\left(\begin{array}{c}
s \\
T\cdot n \\
\end{array}\right)\ dS\]
\[=-\omega\int_{\partial\Gamma}M_2^{-1}\left(\begin{array}{c}
g' \\
g'' \\
\end{array}\right)\cdot\left(\begin{array}{c}
s \\
T\cdot n \\
\end{array}\right)\ dS.\]
Since 
\[M_2^{-1}=\frac{1}{a'}\left(\begin{array}{cc}
1 & 0 \\
-a'' & a' \\
\end{array}\right),\]
we have 
\[M_2^{-1}M_1=\frac{1}{a'}\left(\begin{array}{cc}
1 & -a'' \\
-a'' & |a|^2 \\
\end{array}\right),\]
which is positive definite as long as $a'>0$.  The new bilinear form above is guaranteed to be coercive as long as $\rho$ and $\kappa$ satisfy (\ref{coerc_bounds}) and $a'<0$.

To find a numerical solution for the Robin boundary value problem, we discretize using the finite element scheme presented in Section~\ref{implement}.  Unfortunately, the surface integrals can no longer be converted to volume integrals by integration by parts and must be computed as they stand.  In this case, the finite element matrix is written as the sum of two matrices $A-\omega B$, where $A$ is of the form (\ref{block}), and the blocks have entries (\ref{block_entries}), and another matrix $B$ with the same block form and block entries 
\begin{eqnarray*}(B_1)_{tj}&=&\displaystyle\int_{\partial\Gamma}\left(\begin{array}{c}
 \psi_t \\
0 \\
\end{array}\right)\cdot M_2^{-1}M_1\left(\begin{array}{c}
\psi_j \\
0 \\
\end{array}\right)\ dS \\
(B_2)_{tj}&=&\displaystyle\int_{\partial\Gamma}\left(\begin{array}{c}
0 \\
\phi_{1t}\cdot n \\
\end{array}\right)\cdot M_2^{-1}M_1\left(\begin{array}{c}
0 \\
\phi_{1j}\cdot n \\
\end{array}\right)\ dS \\
(B_3)_{tj}&=&\displaystyle\int_{\partial\Gamma}\left(\begin{array}{c}
0 \\
\phi_{2t}\cdot n \\
\end{array}\right)\cdot M_2^{-1}M_1\left(\begin{array}{c}
0 \\
\phi_{2j}\cdot n \\
\end{array}\right)\ dS \\
(B_4)_{tj}&=&\displaystyle\int_{\partial\Gamma}\left(\begin{array}{c}
0 \\
\phi_{1t}\cdot n \\
\end{array}\right)\cdot M_2^{-1}M_1\left(\begin{array}{c}
\psi_j \\
0 \\
\end{array}\right)\ dS \\
(B_5)_{tj}&=&\displaystyle\int_{\partial\Gamma}\left(\begin{array}{c}
0 \\
\phi_{2t}\cdot n \\
\end{array}\right)\cdot M_2^{-1}M_1\left(\begin{array}{c}
0 \\
\phi_{1j}\cdot n \\
\end{array}\right)\ dS \\
(B_6)_{tj}&=&\displaystyle\int_{\partial\Gamma}\left(\begin{array}{c}
0 \\
\phi_{2t}\cdot n \\
\end{array}\right)\cdot M_2^{-1}M_1\left(\begin{array}{c}
\psi_j \\
0 \\
\end{array}\right)\ dS. \\
\end{eqnarray*}
The right-hand side vector $b$ is also partitioned as $(b_1,b_2,b_3)^T$ with entries
\[\begin{array}{ll}
(b_1)_k= & \displaystyle-\omega\int_{\partial\Gamma}\left(\begin{array}{c}
\psi_k \\
0 \\
\end{array}\right)\cdot M_2^{-1}\left(\begin{array}{c}
g' \\
g'' \\
\end{array}\right)\ dS \\
(b_2)_k= & -\displaystyle\omega\int_{\partial\Gamma}\left(\begin{array}{c}
0 \\
\phi_{1k}\cdot n \\
\end{array}\right)\cdot M_2^{-1}\left(\begin{array}{c}
g' \\
g'' \\
\end{array}\right)\ dS \\
 (b_3)_k= & -\displaystyle\omega\int_{\partial\Gamma}\left(\begin{array}{c}
0 \\
\phi_{2k}\cdot n \\
\end{array}\right)\cdot M_2^{-1}\left(\begin{array}{c}
g' \\
g'' \\
\end{array}\right)\ dS. \\
\end{array}\]

Assuming that the coercivity requirements (\ref{coerc_bounds}) on $\rho$ and $\kappa$ are satisfied, and $a' < 0$, the system 
\[(A-\omega B)\alpha=b\]
may be solved using the same preconditioned conjugate gradient approach as outlined previously.

 \begin{figure}[t]
   \centering
   \includegraphics[width=2.8994in]{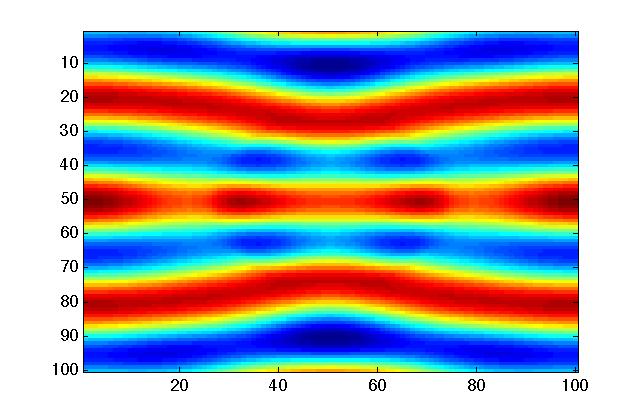} 
   \includegraphics[width=2.8994in]{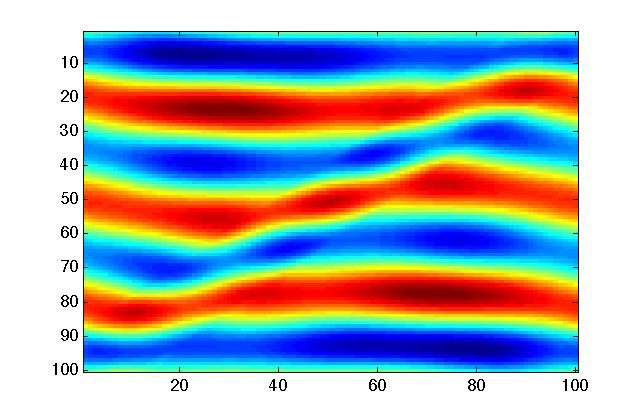} 
   \caption{On the left is shown the solution to the Robin problem with a disc of high density material centered in the domain.  On the right the disk is replaced with a bar of the same material angled from the lower left to the upper right of the domain.}
     \label{circ_scat}
\end{figure}

\subsection{Numerical Examples}

Here we present some numerical examples obtained by using the finite element method to solve problems with Robin boundary conditions.  In these examples the Robin boundary conditions are imposed on $y=0$ and $y=1$, while on the other sides of the domain we have imposed periodic boundary conditions.  On the left in Figure~\ref{circ_scat} is the  solution with a circular scatterer with $\rho=1+.011i$ outside the scatterer, $\rho=2+.011i$ inside the scatterer, $\kappa=1+.011i$ everywhere, $a=-1+.333i$ and $g=3.33i$.  On the right, the circular scatterer is replaced by a bar angled across the domain, but the other parameters in the problem remain the same.  These results were calculated using the $RT_{[0]}$ discretization for the $v$ variable described in Section~\ref{other_disc}.

\section{Conclusions}
The variational principles of Milton, Seppecher, Bouchitt\'{e}, and Willis make it possible to formulate the solution of the Helmholtz equation as a minimization, and this is reflected in the fact that the stiffness matrix for the finite element method is symmetric positive definite.  This allows us to use classical iterative methods such as preconditioned conjugate gradient to solve the associated system, along with straightforward finite element error estimates.  The primary advantage of this approach is that it allows the use of efficient iterative methods for the solution of the linear system. But there are also disadvantages in that the system has more unknowns, since we must solve for $P$ and $v$ simultaneously.  

More research is necessary to determine the circumstances under which this approach may be more effective than others currently in use.  A particular point of interest is that even though the underlying minimization principles are valid for arbitrarily small loss coefficients, the conditioning of the associated finite element matrix deteriorates as the system becomes less dissipative.  The general question of how solution efficiency depends on loss should be analyzed further.

We have only approached the scalar, two-dimensional Helmholtz equation in this paper, while the minimization principles of Milton, Seppecher, Bouchitt\'{e}, and Willis apply to the full vector Maxwell equations, as well as the equations of linear elasticity.  The general approach taken here should also apply in those cases.  We note finally that in many applications for these models, one would like to apply nonlocal transmission or radiation boundary conditions in order to accurately handle unbounded domains.  The problem of adapting these minimization methods to such boundary conditions remains open, although presumably the PML approach (see eg. \cite{PML}) would apply.

\section*{Acknowledgements}

Russell Richins is grateful for support from the National Science Foundation through grant DMS-0707978.  Also, the authors would like to thank Graeme Milton for many helpful suggestions, and John Willis who, along with Graeme Milton, clarified an earlier formulation of the boundary value problem for us.

\bibliographystyle{amsplain}
\bibliography{paperbibliography}

\providecommand{\bysame}{\leavevmode\hbox to3em{\hrulefill}\thinspace}
\providecommand{\MR}{\relax\ifhmode\unskip\space\fi MR }
\providecommand{\MRhref}[2]{%
  \href{http://www.ams.org/mathscinet-getitem?mr=#1}{#2}
}
\providecommand{\href}[2]{#2}
\begin{thebibliography}{1}

\bibitem{Axelsson-Barker}
O.~Axelsson and V.~A. Barker, \emph{Finite element solution of boundary value
  problems, theory and computation}, SIAM, Philidelphia, PA, 2001.

\bibitem{Brezzi-Fortin}
F.~Brezzi and M.~Fortin, \emph{Mixed and hybrid finite element methods},
  Springer-Verlag, New York, NY, 1991.

\bibitem{Cherk}
A.V. Cherkaev and L.V. Gibiansky, \emph{Variational principles for complex
  conductivity, viscoelasticity, and similar problems in media with complex
  moduli}, Journal of Mathematical Physics \textbf{35} (1994), 127--145.

\bibitem{Demmel}
J.~Demmel, \emph{The condition number of equivalence transformations that block
  diagonalize matrix pencils}, S{I}{A}{M} J. Num. Anal. \textbf{20} (1983),
  599--610.

\bibitem{Demmel_book}
\bysame, \emph{Applied numerical linear algebra}, {S}{I}{A}{M}, Philidelphia,
  PA, 1997.

\bibitem{Evans}
L.C. Evans, \emph{Partial differential equations}, American Mathematical
  Society, Providence, RI, 1998.

\bibitem{Graeme}
G.W. Milton, P.~Seppecher, and G.~Bouchitt\'{e}, \emph{Minimization variational
  principles for acoustics, elastodynamics, and electromagnetism in lossy
  inhomogeneous bodies at fixed frequency}, Proc. R. Soc. A \textbf{465}
  (2009), 367--396.

\bibitem{Willis}
G.W. Milton and J.R. Willis, \emph{Minimum variational principles for
  time-harmonic waves in a dissipative medium and associated variational
  principles of {H}ashin-{S}htrikman type}, Proc. Roy. Soc. Lond. to appear.

\bibitem{PML}
I.~Harari~M. Slavutin and E.~Turkel, \emph{Analytical and numerical studies of
  a finite element {PML} for the {H}elmholtz equation}, J. Comp. Acoust.
  \textbf{8} (2000), 121--137.

\end{thebibliography}

\end{document}